\documentclass{article}

\newtheorem{theorem}{Theorem}[section]
\newtheorem{definition}{Definition}

\begin{document}

\title{The Poincar\'e conjecture \\  for stellar manifolds}

\author{Sergey Nikitin \\
Department of Mathematics\\
Arizona State University\\
Tempe, AZ 85287-1804\\
 nikitin@asu.edu \\
 http://lagrange.la.asu.edu}

\maketitle

{\bf Abstract}  This paper proves that any closed, simply connected, connected, compact stellar manifold is a stellar sphere. That implies the Poincar\'e conjecture.

\vspace{0.1cm}

\section{Introduction}

\vspace{0.1cm}

We prove that a closed, compact, connected and simply connected stellar manifold is a stellar sphere. As a corollary we obtain the
famous Poincar\'e conjecture:

\vspace{0.1cm}

{\it
Every simply connected closed $3$-manifold is homeomorphic to the $3$-sphere.
} 

\vspace{0.1cm}
That was stated by Henri Poincar\'e in 1904 \cite{Poincare}.
Analogues of this hypothesis were successfully proved in dimensions higher than $3,$
see  \cite{Freedman}, \cite{Milnor}, \cite{Smale},\cite{Wallace}, \cite{Zeeman}.
\\
We prove this conjecture for stellar manifolds.
Since every $3$-dimensional manifold can be triangulated and any two stellar equivalent manifolds are piecewise linearly homeomorphic (\cite{Lickorish},\cite{Newman},\cite{Pachner}), our result does imply the famous Poincar\'e conjecture.

\section{Main result}
A stellar $n$-manifold $M$ can be identified with the sum of its $n$-dimensional simplexes ($n$-simplexes):
$$
M=\sum_{i=1}^n g_i
$$
with coefficients from ${\rm Z}_2.$ We will call $\{g_i\}_{i=1}^n$ generators of $M.$

All vertices in $M$ can be enumerated and any $n$-simplex $s$ from $M$ corresponds to the set
of its vertices
$$
s=(i_1 \; i_2 \; \dots \; i_{n+1}),
$$
where $i_1 \; i_2 \; \dots \; i_{n+1}$ are integers. 

The boundary operator $\partial $ is defined on a simplex as
$$
\partial (i_1 \; i_2 \; \dots \; i_{n+1}) = (i_1 \; i_2 \; \dots \; i_n ) + (i_1 \; i_2 \; \dots \; i_{n-1} \; i_{n+1} ) + \dots + (i_2 \; i_3 \; \dots \; i_{n+1})
$$
and linearly extended to any complex, i.e.
$$
\partial M = \sum_{i=1}^n \partial g_i.
$$
A manifold is called closed if $\partial M =0.$

If two simplexes $(i_1 \; i_2 \; \dots i_m) $ and $(j_1 \; j_2 \; \dots j_n)$ do not have common vertices then one can define their join
$$
(i_1 \; i_2 \; \dots i_m) \star (j_1 \; j_2 \; \dots j_n)
$$
as the union
$$
(i_1 \; i_2 \; \dots i_m) \cup (j_1 \; j_2 \; \dots j_n).
$$

\vspace{0.1cm}

If two complexes $K=\sum_i q_i $ and $ L = \sum_j p_j$ do not have common vertices then 
their join is defined as
$$
K\star L = \sum_{i,j} q_i \star p_j.
$$

If $A$ is a simplex in a complex $K$ then we can introduce its link:
$$
lk(A,K) = \{ B \in K \; ; \; A \star B \in K \}.
$$
The star of $A$ in $K$ is $A \star lk(A,K).$ Thus,
$$
K = A\star lk(A,K) + Q(A,K),
$$
where the complex $Q(A,K)$ is composed of all the generators of $K$ that do not contain $A.$ A complex with generators of the same dimension is called a uniform complex.

\begin{definition} ({\bf Subdivision})
Let $A$ be a simplex of a complex $K.$ Then any integer $a$ which is not a vertex of $K$ defines starring of
$$
K=A\star lk(A,K) + Q(A,K)
$$
at $a$ as
$$
 \hat{K}=a\star \partial A \star lk(A,K) + Q(A,K).
$$
This is denoted as
$$
\hat{K} =  (A \; a)K.
$$
\end{definition}

The next operation is the inverse of subdivision. It is called a stellar weld and defined as follows.

\begin{definition} ({\bf Weld})
  Consider a complex
  $$
        \hat{K}=a\star  lk(a,\hat{K}) + Q(a,\hat{K}),
  $$
with $lk(a,\hat{K}) = \partial A \star B$ where
$B$ is a subcomplex in $\hat{K},$  $A$ is a simplex and $A\notin \hat{K} .$
Then the (stellar) weld $(A\; a)^{-1} \hat{K}$ is defined as
  $$
       (A\; a)^{-1} \hat{K} =  A \star B  + Q(a,\hat{K}).
  $$
\end{definition}

A stellar move is one of the following operations: subdivision, weld, enumeration change on the set of vertices.
Two complexes $M$ and $L$  are called stellar equivalent if one is obtained from the other by a finite sequence of stellar moves.
It is denoted as $M \sim L.$\\

If a complex $L$ is stellar equivalent to $(1 \; 2 \; \dots \; n+1)$ then $L$ is called a stellar 
$n$-ball. On the other hand, if $K \sim \partial (1 \; 2 \; \dots \; n+2)$ then $K$ is a stellar
$n$-sphere.

\begin{definition} ({\bf Stellar manifold})
\label{stellar_def}
Let $M$ be a complex. If, for every vertex $i$ of $M,$ the link $lk(i,M)$ is either a stellar $(n-1)$-ball or a stellar $(n-1)$-sphere, then $M$ is  a
stellar $n$-dimensional manifold ($n$-manifold).
\end{definition}

If $i$ is a vertex of $M$ then
$$
M= i\star lk(i,M) + Q(i,M).
$$
If $\partial M=0,$ then $Q(i,M)$ is a stellar manifold.\\
Indeed, consider an arbitrary vertex $j$ of $Q(i,M).$ Then
$$
lk(i,M)=j\star lk(j,lk(i,M)) + Q(j,lk(i,M))
$$
and
$$
Q(i,M)=j\star lk(j,Q(i,M)) + Q(j,Q(i,M)).
$$
Since $M$ is a stellar manifold and $\partial M = 0$
$$
i\star lk(j,lk(i,M)) + lk(j,Q(i,M)) 
$$
is a stellar sphere. Hence, it follows from \cite{Newman} that $lk(j,Q(i,M))$ is either a stellar ball or a stellar sphere.\\
In the sequel it is convenient to consider an equivalence relation on the set of vertices of a stellar manifold. Among all possible such 
equivalence relations we are mostly interested in those that meet certain regularity properties underlined by the following definition.
\begin{definition} ({\bf Regular equivalence})
Given a stellar manifold $M,$ an equivalence relation "$\simeq$" on the set of vertices from $M$ is called regular if it meets the following conditions:
\begin{itemize}
\item[(i)] No generator $g \in M$ has two vertices that are equivalent to each other.
\item[(ii)]For any generator $g\in M$ there might exist not more than one generator $p\in M \setminus g$ such that any vertex of $g$ is either equal or equivalent to some vertex of $p.$ We call such two generators equivalent, $g \simeq p.$ 
\end{itemize}
\end{definition}
Our proof of the Poincar\'e conjecture is based on the following result.

\begin{theorem}
\label{triangulation}
A connected stellar $3$-manifold $M$ with a finite number of generators admits a triangulation
$$
N = a\star (S/\simeq) ,
$$
where $a\notin S$ is a vertex, $S$ is a stellar $2$-sphere and "$\simeq $" is a regular equivalence relation. Moreover, if $M$ is closed then for any generator $g\in S$ there exists
exactly one generator $p\in S\setminus g$ such that $g \simeq p.$
\end{theorem}
{\bf Proof.}
Let us choose an arbitrary generator $g \in M$ and an integer $a$ that is not a vertex of $M.$ Then
$$
M \sim (g \; a ) M  \mbox{  and  }
(g \; a ) M =a \star \partial g + M\setminus g,
$$
where $M\setminus g$ is defined by all the generators of $M$ excluding $g.$
We construct $N$ in finite number of steps.
Let $N_0=(g \; a ) M.$ Suppose we constructed already $N_k$ and there exists a generator $p\in  Q(a,N_k)$ that has at least one common $2$-simplex with $lk(a, N_k).$
Without loss of generality, we can assume that 
$$
 p = (1 \; 2\; 3 \; 4).
$$
and $(1 \; 2\; 3)$ belongs to $lk(a, N_k).$
If the vertex $(4)$ does not belong to  $lk(a, N_k)$ then
$$
N_{k+1} = ((a\; 4 )\; b)^{-1}((1 \; 2 \; 3) \; b) N_k.
$$
If the vertex $(4)$ belongs to $lk(a, N_k)$ then after introducing a new vertex $d\notin N_k$ we take
$$
L = ((a\; d )\; b)^{-1}((1 \; 2 \; 3) \; b)(N_k\setminus p + (1\;2\;3\;d)),
$$
where $b\notin (N_k\setminus p + (1\;2\;3\;d)),$ and
$$
N_{k+1} =  a \star (lk(a,L)/\simeq ) + Q(a,L)
$$
endowed with the equivalence $d \simeq (4).$\\
By construction
$$
N_{k+1} = a \star (lk(a,N_k) +\partial p ) + Q(a,N_k) \setminus p
$$
if $(4)\notin lk(a, N_k).$ Otherwise,
$$
N_{k+1} = a \star ((lk(a,N_k) +\partial g)/\simeq ) + Q(a,N_k) \setminus p,
$$
where $g=(1 \; 2\; 3 \; d),\;\;d\simeq (4).$

Since $M$ is connected and has a finite number of generators there exists a natural number $m$ such that 
$$
N_m =  a \star (S/\simeq)
$$
where $S$ is a stellar $2$-sphere and "$\simeq $" is a regular equivalence relation.\\
If $M$ is closed, then $\partial N_m = 0,$ and therefore, for any generator $g \in S$ there exists exactly one generator $p\in S \setminus g$ such that $g \simeq p.$\\
Q.E.D.\\

It is known \cite{Glaser}, \cite{Lickorish}, \cite{Moise1} that any two $3$-dimensional manifolds admitting stellar equivalent triangulations
are piecewise linearly homeomorphic. On the other hand, every compact $3$-dimensional manifold admits a stellar triangulation with
a finite number of generators \cite{Moise1}. Hence, the Poincar\'e conjecture
 follows from the following statement.

\begin{theorem} ({\bf the Poincar\'e Conjecture})
A simply connected, connected closed stellar $3$-manifold $M$ with a finite number of generators is homeomorphic to the $3$-sphere, $\partial (1\;2\;3\;4\;5).$
\end{theorem}
{\bf Proof.}
Let $M$ be a closed, connected and simply connected $3$-dimensional stellar manifold with a finite number of generators. By Theorem \ref{triangulation} $M$ admits a triangulation
$$
N = a\star (S/\simeq ) 
$$
where $a \notin S $ is a vertex, $S$ is a stellar $2$-sphere and  "$\simeq $" is a regular equivalence relation. Moreover, for any generator $g\in S$ there exists
exactly one generator $p\in S\setminus g $ such that $g\simeq p.$\\
Let us show that $a\star (S/\simeq )$ is homeomorphic to the $3$-sphere.
If $g$ is a generator of $S,$ then 
$$
a \star (S\setminus g)/\simeq
$$
is connected and, by Seifert -- Van Kampen theorem (see e.g. \cite{Massey}), it is simply connected.\\
The barycentric subdivision $Br(a \star (S\setminus g)/\simeq)$ of $a \star (S\setminus g)/\simeq$ is geometrically collapsible (or $Br(a \star (S\setminus g)/\simeq)\searrow 0$). Thus by Whitehead's result on regular neighborhoods (see \cite{Glaser}, \cite{Whitehead})  $Br(a \star (S\setminus g)/\simeq)$ is a combinatorial $3$-ball. Therefore $N = a\star (S/\simeq ),$ two $3$-balls identified  along their boundaries, and $N$ is homeomorphic to the $3$-sphere. Thus, $M$ is homeomorphic to the $3$-sphere. 
Q.E.D.\\

\section{Acknowledgement}
The author is grateful to Prof. C. Margerin and Prof. W.B.Raymond Lickorish for reading early versions of this manuscript and making valuable suggestions that helped to correct mistakes and improve the presentation of the main result of this paper.

\bibliographystyle{amsplain}

\end{document}